\documentclass[12pt]{amsart}
\usepackage{graphicx}
\usepackage{amssymb}
\usepackage{amsfonts}

\usepackage{tikz}
\usetikzlibrary{positioning,arrows,matrix,chains,shadows,shapes}
\tikzstyle{Point} = [fill, radius=0.08]

\usepackage{pspicture}
\usepackage{pstricks}
\usepackage[matrix,frame,xdvi]{xypic} 

\newdimen\vcadre\vcadre=0.1cm 
\newdimen\hcadre\hcadre=0.1cm 
\def\GrTeXBox#1{\vbox{\vskip\vcadre\hbox{\hskip\hcadre%
      $#1$%
   \hskip\hcadre}\vskip\vcadre}}
\def\arx#1[#2]{\ifcase#1 \relax \or%
  \ar @{-}[#2]  \or%
  \ar @2{-}[#2] \or%
  \ar @{--}[#2] \or%
  \ar @2{.}[#2] \or%
  \ar @{~}[#2]  \fi}

\def\arbg#1#2{
\newdimen\vcadre\vcadre=0.01cm 
\newdimen\hcadre\hcadre=0.01cm 
\xymatrix@R=0.1cm@C=2mm{
 && {\GrTeXBox{#1}}\arx1[ld]\arx1[rd]\\
& {\GrTeXBox{#2}}\arx1[dl]\arx1[dr] && {\GrTeXBox{z}}\\
 {\GrTeXBox{x}} &&  {\GrTeXBox{y}}\\
}
}
\def\arbd#1#2{
\newdimen\vcadre\vcadre=0.01cm 
\newdimen\hcadre\hcadre=0.01cm 
\xymatrix@R=0.1cm@C=2mm{
 & {\GrTeXBox{#1}}\arx1[ld]\arx1[rd]\\
{\GrTeXBox{x}} & & {\GrTeXBox{#2}}\arx1[ld]\arx1[rd]\\
& {\GrTeXBox{y}}&  & {\GrTeXBox{z}}\\
}
}
\def\arbgnew#1#2{
\newdimen\vcadre\vcadre=0.01cm 
\newdimen\hcadre\hcadre=0.01cm 
\xymatrix@R=0.1cm@C=2mm{
 && {\GrTeXBox{#1}}\arx1[ld]\arx1[rd]\\
& {\GrTeXBox{#2}}\arx1[dl]\arx1[dr] && {\GrTeXBox{B}}\\
 {\GrTeXBox{.}} &&  {\GrTeXBox{A}}\\
}
}

\newtheorem{example}{Example}[section]
\newtheorem{note}[example]{Note}
\newtheorem{theorem}[example]{Theorem}

\newtheorem{proposition}[example]{Proposition}

\def\Proof{\noindent \it Proof -- \rm}
\def\qed{\hspace{3.5mm} \hfill \vbox{\hrule height 3pt depth 2 pt width 2mm}
\bigskip}

\def\K{{\mathbb K}}

\def\PBT{{\bf PBT}}

\def\FQSym{{\bf FQSym}}

\def\WQSym{{\bf WQSym}}
\def\ev{{\rm ev}}

\def\ssh{\Cup}

\def\sconc{\bullet}

\def\<{\langle}
\def\>{\rangle}

\def\F{{\bf F}}

\def\dim{{\rm dim}}

\def\mPBT{{}^m{\bf PBT}}

\def\mFQSym{{}^m{\bf FQSym}}

\def\D{{\mathcal D}}

\def\shuff#1#2{\mathbin{
\hbox{\vbox{ \hbox{\vrule \hskip#2 \vrule height#1 width 0pt
}%
\hrule}%
\vbox{ \hbox{\vrule \hskip#2 \vrule height#1 width 0pt
\vrule }%
\hrule}%
}}}

\def\binomial#1#2{\left(\,\begin{matrix}#1 \cr #2\end{matrix}\,\right)}

\def\shuf{{\mathchoice{\shuff{7pt}{3.5pt}}%
{\shuff{6pt}{3pt}}%
{\shuff{4pt}{2pt}}%
{\shuff{3pt}{1.5pt}}}}%
\def\shuffle{\,\shuf\,}

\def\ddeux{D^{(2)}}
\def\DP#1{Dp^{(#1)}}

\def\gaudend{\prec}     
\def\droitdend{\succ}   

\def\gd{\prec}
\def\md{\circ}
\def\dd{\succ}

\def\dgd{\dashv}   
\def\dmd{\perp} 
\def\ddd{\vdash}   

\def\bulle{\bullet}

\tikzset{baseline={([yshift=-3.5pt]current bounding box.center)}}
\tikzset{level distance = 0.7cm, sibling distance = 2em}
\tikzset{edge from parent/.style={
       draw,
       edge from parent path = {(\tikzparentnode) --
                                (\tikzchildnode)}
    }}
\tikzset{root/.style={}}
\tikzset{white/.style={}}

\def\treeLft#1#2{\begin{tikzpicture}
\node [root] {$#1$}
  child {
    node [white] {$#2$}
      child [missing]
      child [missing]
  }
  child [missing];
\end{tikzpicture}}

\def\treeRgt#1#2{\begin{tikzpicture}
\node [root]{$#1$}
  child [missing]
  child {
    node [white] {$#2$}
      child [missing]
      child [missing]
  };
\end{tikzpicture}}

\def\treeGen#1#2#3{\begin{tikzpicture}
\node [root]{$#1$}
  child {
    node [white] {$#2$}
      child [missing]
      child [missing]
  }
  child {
    node [white] {$#3$}
      child [missing]
      child [missing]
  };
\end{tikzpicture}}

\def\treeGenG#1#2#3#4#5{\begin{tikzpicture}
\node [root]{$#1$}
  child {
    node [white] {$#2$}
      child {node [white] {$#3$} child [missing] child [missing]}
      child {node [white] {$#4$} child [missing] child [missing]}
  }
  child {
    node [white] {$#5$} child [missing] child [missing]
  };
\end{tikzpicture}}

\def\treeGenD#1#2#3#4#5{\begin{tikzpicture}
\node [root]{$#1$}
  child {
    node [white] {$#2$} child [missing] child [missing]
  }
  child {
    node [white] {$#3$}
      child {node [white] {$#4$} child [missing] child [missing]}
      child {node [white] {$#5$} child [missing] child [missing]}
  };
\end{tikzpicture}}

\def\treeGenT#1#2#3#4{\begin{tikzpicture}
\node [root]{$#1$}
  child { node [white] {$#2$} }
  child { node [white] {$#3$} }
  child { node [white] {$#4$} };
\end{tikzpicture}}

\def\treeGenTd#1#2#3{\begin{tikzpicture}
\node [root]{$#1$}
  child { node [white] {$#2$} }
  child { node [white] {$#3$} }
  child [missing];
\end{tikzpicture}}

\def\treeGenTz#1#2{\begin{tikzpicture}
\node [root]{$#1$}
  child { node [white] {$#2$} }
  child [missing]
  child [missing];
\end{tikzpicture}}

\def\treeGenTs#1#2#3#4#5#6{\begin{tikzpicture}
\node [root]{$#1$}
  child { node [white] {$#2$} }
  child { node [white] {$#3$}
         child { node [white] {$#4$} }
         child { node [white] {$#5$} }
         child { node [white] {$#6$} }
        }
  child [missing];
\end{tikzpicture}}

\def\treeGenTes#1#2#3#4#5#6#7{\begin{tikzpicture}
\node [root]{$#1$}
  child { node [white] {$#2$} }
  child { node [white] {$#3$} }
  child { node [white] {$#4$}
         child { node [white] {$#5$} }
         child { node [white] {$#6$} }
         child { node [white] {$#7$} }
        };
\end{tikzpicture}}

\def\treeEx{\begin{tikzpicture}
\node [root]{$\bulle$}
  child { node [white] {$\bulle$} }
  child { node [white] {$\bulle$}
    child { node [white] {$\bulle$} 
      child [missing]
      child { node [white] {$\bulle$}
      child [missing]
      child [missing]
      child { node [white] {$\bulle$} }
      }
      child [missing]
    }
    child { node [white] {$\bulle$} }
    child { node [white] {$\bulle$} }
  }
  child { node [white] {$\bulle$} };
\end{tikzpicture}}
\def\treeExlab{\begin{tikzpicture}
\node [root]{$8$}
  child { node [white] {$7$} }
  child { node [white] {$5$}
    child { node [white] {$4$}
      child [missing]
      child { node [white] {$2$}
      child [missing]
      child [missing]
      child { node [white] {$3$} }
      }
      child [missing]
    }
    child { node [white] {$1$} }
    child { node [white] {$6$} }
  }
  child { node [white] {$9$} };
\end{tikzpicture}}

\def\treeExb{\begin{tikzpicture}
\node [root]{$\bulle$}
  child { node [white] {$\bulle$} }
  child { node [white] {$\bulle$}
    child { node [white] {$\bulle$}
      child [missing]
      child { node [white] {$\bulle$}
      child [missing]
      child [missing]
      child { node [white] {$\bulle$} }
      }
      child [missing]
    }     
    child { node [white] {$\bulle$} }
    child { node [white] {$\bulle$} }
  }     
  child [missing];
\end{tikzpicture}}

\def\treeExc{\begin{tikzpicture}
\node [root]{$\bulle$}
    child { node [white] {$\bulle$}
      child [missing]
      child { node [white] {$\bulle$}
      child [missing]
      child [missing]
      child { node [white] {$\bulle$} }
      }
      child [missing]
    }
    child { node [white] {$\bulle$}}
  child [missing];
\end{tikzpicture}}
\def\treeExd{\treeGenTd\bulle\bulle\bulle}
\def\treeExe{\begin{tikzpicture}
\node [root]{$\bulle$}
    child { node [white] {$\bulle$}
      child [missing]
      child { node [white] {$\bulle$}
      child [missing]
      child [missing]
      child { node [white] {$\bulle$} }
      }
      child [missing]
    }
    child [missing]
  child [missing];
\end{tikzpicture}}
\def\treeExf{\begin{tikzpicture}
\node [root]{$\bulle$}
      child [missing]
      child { node [white] {$\bulle$}
      child [missing]
      child [missing]
      child { node [white] {$\bulle$} }
      }
  child [missing];
\end{tikzpicture}}
\def\treeExg{\begin{tikzpicture}
\node [root]{$\bulle$}
  child [missing]
  child { node [white] {$\bulle$}}
  child [missing];
\end{tikzpicture}}
\def\treeExh{\begin{tikzpicture}
\node [root]{$\bulle$}
  child { node [white] {$\bulle$}}
  child [missing]
  child [missing];
\end{tikzpicture}}

\begin{document}
\title[$m$-dendriform algebras]%
{$m$-dendriform algebras}
\author[J.-C.~Novelli]%
{Jean-Christophe Novelli}

\address[]{\hskip-.5cm Laboratoire d'informatique Gaspard-Monge\\
Universit\'e Paris-Est Marne-la-Vall\'ee \\
5, Boulevard Descartes \\ Champs-sur-Marne \\
77454 Marne-la-Vall\'ee cedex 2 \\
France}
\email[Jean-Christophe Novelli]{novelli@univ-mlv.fr}
\date{\today}


\begin{abstract}
The Fuss-Catalan numbers are a generalization of the Catalan numbers.
They enumerate a large class of objects and in particular $m$-Dyck paths and
$(m\!+\!1)$-ary trees.
Recently, F. Bergeron defined an analogue for generic $m$ of the Tamari order
on classical Dyck words. The author and J.-Y. Thibon showed that the
combinatorial Hopf algebras related to these $m$-Tamari orders are defined
thanks to the same monoid, the sylvester monoid, as in the $m=1$ case and that
all related Hopf algebras also have $m$ analogues.

We present here the $m$-generalization of another construction on Catalan
sets: the dendriform algebras. These algebras are presented in two different
ways: first by relations between the $m\!+\!1$ operations, relations that are
very similar to the classical relations; and then by explicit operations
splitting the classical dendriform operations defined on words into new
operations. We then investigate their dual and show they are Koszul.
\end{abstract}

\maketitle
\footnotesize
\tableofcontents
\normalsize

\section{Introduction}

Among the so-called combinatorial Hopf algebras, the Loday-Ronco algebra of
planar binary trees \cite{LR1} plays a prominent role.
Its original definition comes from uhe theory of operads, as this algebra
is the free dendriform algebra on one generator. One can also define it in a
combinatorial way by means of the so-called sylvester congruence~\cite{HNT}
applied to the Malvenuto-Reutenauer algebra $MR$~\cite{MR}, an algebra defined
on permutations.
Since the product in this algebra consists in shuffling permutations, it is
immediate to check that all products are intervals of the weak order on
permutations. One very nice property of the Loday-Ronco algebra is that
the congruence classes also are intervals of the weak order, so that products
in this algebra also are intervals of the weak order. One can then consider a
suborder by restricting to particular class representatives. As Loday and
Ronco pointed out, this suborder appears to be a well-known order, the Tamari
order initially defined as an order of parenthesized expressions.

\medskip
Recently, the Tamari order has been generalized by Bergeron to an infinite
series of lattices, the $m$-Tamari orders \cite{Ber}, defined on combinatorial
objects counted by Fuss-Catalan numbers, for example $(m+1)$-ary trees. This
raises the question of the existence of a generalization of the dendriform
algebra for all values of $m$.

We shall give a positive answer to this question. Actually, in~\cite{NT14},
the authors showed that there is a way to generalize permutations so
that the number of sylvester classes of these $m$-permutations is the same
as the number of $(m+1)$-ary trees. So, a way to look for an $m$-dendriform
algebra is to start with the algebra on $m$-permutations and split the product
into $m+1$ operations that preserve sylvester classes and such that all
sylvester classes can be obtained by a suitable sequence of operations
applied to one generator. We shall see that there is a very simple way to
define these operations using the right-to-left minima on words.
Note that Leroux also defined algebras that deserve the name of
$m$-dendriform algebras in~\cite{Ler} very similar to ours but not isomorphic
and of a different spirit: the author gave the relations without a
realization. More details can be found in Note~\ref{note-leroux}.

Finally, the initial study of dendriform algebras was deeply connected to the
study of their dual, the $Dias$ algebras that are algebras satisfying monoidal
relations between their operations. We shall compute here the dual of the
$m$-dendriform algebras, hence generalizations of the $Dias$ algebras.

\bigskip
This article is structured as follows. We first fix our notations on words and
their operations and recall the basic definitions and properties of the
sylvester congruence. We then recall some facts about dendriform algebras and
sketch our way to prove those. We then present in full detail the case $m=2$
of the construction and show that our $2$-dendriform algebra has as series of
dimensions the generating series of ternary trees.
We then present the case for general $m$ which almost consists in copying the
$m=2$ case. We end by presenting the structure of the dual algebras of our
$m$-dendriform algebras and show that their operads are Koszul operads.

\bigskip
{\bf Acknowledgements.} Partially supported by ANR CARMA. The author also
wishes to thank Samuele Giraudo for fruitful discussions about Koszul duality.

\section{Background and notations}

\subsection{Words}

In the sequel, we shall need a countable totally ordered alphabet $A$, usually
labeled by the positive integers.
We denote by $A^*$ the free monoid generated by $A$.

All algebras will be taken over a field $\K$ of characteristic $0$. The
notation $\K\<A\>$ means the free associative algebra over $A$ when $A$ is
finite, and the inverse limit $\varprojlim_B \K\<B\>$, where $B$ runs over
finite subsets of $A$, when $A$ is infinite.

The \emph{evaluation} $\ev(w)$ of a word $w$ is the sequence whose $i$-th term
is the number of occurrences of the letter $a_i$ in $w$.

\medskip
For a word $w$ over the alphabet $\{1,2,\ldots\}$, we denote by $w[k]$ the
word obtained by replacing each letter $i$ by the integer $i+k$.
If $u$ and $v$ are two words, with $u$ of length $k$, one defines
the \emph{shifted concatenation}
$u\sconc v = u\cdot (v[k])$
and the \emph{shifted shuffle}
$ u\ssh v= u\shuffle (v[k])$,
where $w_1\shuffle w_2$ is the usual shuffle product defined recursively by
\begin{itemize}
\item $w_1 \shuffle \epsilon  = w_1$,\quad $\epsilon \shuffle w_2  = w_2$,
\item $au \shuffle bv = a(u \shuffle bv) + b(au \shuffle v)$,
\end{itemize}
where $w_1=a\cdot u$ and $w_2=b\cdot v$, and both $a$ and $b$ are letters and
$\cdot$ means concatenation.

For example,
\begin{equation}
12\ssh 21 = 12\shuffle 43 = 1243 + 1423 + 1432 + 4123 + 4132 + 4312\,.
\end{equation}

\subsection{The sylvester monoid}

\subsubsection{Definition}

The sylvester monoid has been defined in~\cite{HNT} to show the parallel
between the algebra of planar binary trees defined by Loday and Ronco
in~\cite{LR1} and the algebra of standard Young tableaux defined
in~\cite{NCSF6}.
Indeed, the sylvester monoid has essentially the same properties as the
plactic monoid and both presentations are indeed very close.

Let us recall that two words $w_1$ and $w_2$ are sylvester-adjacent if there
exists three words $u$, $v$, $w$ and three letters
$a\leq b<c$ such that
\begin{equation}
w_1 = u\ ac\ v\ b\ w \text{\qquad and\qquad}
w_2 = u\ ca\ v\ b\ w.
\end{equation}
The \emph{sylvester congruence} is then the transitive closure of the
sylvester adjacence and the \emph{sylvester monoid} is the quotient of the
free monoid $A^*$ by the sylvester congruence, where $A$ is any ordered
alphabet.

There are many well-known facts about the sylvester monoid and we shall only
recall here the ones we shall use. The reader can find more details and
more results in~\cite{HNT,NT14}.

First, any sylvester class is in bijection with a binary search tree, the
elements of the class then being its linear extensions. Moreover, the
sylvester classes of permutations form an interval of the weak order on
permutations, and the greatest elements of the sylvester classes are the
permutations avoiding the pattern $132$.

In particular, the sylvester classes of permutations of size $n$ are in
bijection with unlabelled binary trees of size $n$ since given such a tree,
there is a unique labelling using each integer from $1$ to $n$ once,
that provides a binary search tree. So sylvester classes of permutations
of size $n$ are enumerated by the Catalan numbers and are more precisely
related to the well-known Tamari order since it is the order obtained by
selecting in the weak order the permutations avoiding the pattern $132$.

Since the Tamari order was generalized by Bergeron~\cite{Ber} where binary
trees are replaced by $(m+1)$-ary trees, it was tempting to look for another
set of words whose sylvester classes are enumerated by the number of such
trees, which are the Fuss-Catalan numbers.
It happens that one can enumerate in an efficient way the number of sylvester
classes of all words of a given evaluation~\cite{NT14} but the final formula
is not as simple as the Catalan numbers, except on the case we shall now
detail.

\subsubsection{Sylvester classes of $m$-permutations and $(m\!+\!1)$-ary trees}
\label{sec-sylvmp}

Indeed, there is one more general case where there is a nice formula for the
number of sylvester classes of a given evaluation, that is, for evaluations of
the form $m^n$. We shall call these elements the $m$-permutations since each
integer from $1$ to $n$ appears exactly $m$ times.

In that case, the binary search trees of this evaluation correspond
to special binary trees called $m$-binary trees in~\cite{CP} since they hide
inside their structure a recursive $m$-structure. Ch\^atel and Pons provided a
bijection between these trees and the $(m+1)$-ary trees hence showing in a
direct combinatorial way that the number of sylvester classes of evaluation
$m^n$ is the Fuss-Catalan number
\begin{equation}
C^{(m)}_n = \frac{1}{mn+1}\binom{mn+n}{n}.
\end{equation}
Note that $m=1$ gives back the Catalan numbers.

\begin{figure}[ht]
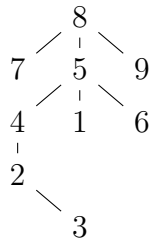

\treeExlab
\caption{\label{fig-tern}The unique labelling of a ternary tree.}
\end{figure}

We shall not need the $m$-binary trees but we shall at some point make use of
ternary trees and even $(m+1)$-ary trees.
Let us define a \emph{$(m+1)$-ary search tree} as the filling of such a tree
by the integers from $1$ to $n$ each once in the following recursive way:
label in increasing order first the $m$-th subtree, then the $(m\!-\!1)$-st
subtree, and so on up to the first subtree, then the root and finally the
$(m\!+\!1)$-st subtree.
For example, Figure~\ref{fig-tern} presents the labelling of a ternary tree
($m=2$ case).

Now, the bijection between the $(m+1)$-ary search trees and the
$m$-permutations avoiding the pattern $132$ is the following: denote by $r$
the root of the tree and read the tree recursively by first reading its right
subtree, then its left subtree, then put $r$, then the second subtree, then
put $r$, and so on up to the $m$-th subtree and the root one last time.

For example, the tree represented Figure~\ref{fig-tern} corresponds to the
word
\begin{equation}
9977 8 66433224 5115 8.
\end{equation}
One easily rebuilds the tree from the word since its root is labelled by the
last letter, which gives back non ambiguously its subtrees since they are all
separated by the other occurrences of the last letter except for the left and
the right subtrees that can be dintinguished by the fact that the letters of
the right subtree are the letters greater than the root.

\subsection{The Hopf algebra $\FQSym$ and its generalizations}

As Loday and Ronco made use in~\cite{LR1} of the algebra of permutations
defined by Malvenuto and Reutenauer, we shall make use of other algebras
defined in~\cite{NT14}. Recall that the Malvenuto-Reutenauer algebra was
originally defined in an abstract way~\cite{MR} and was later proved to be
isomorphic to the algebra whose basis elements $\F_\sigma$ are given by sums
of words having as result $\sigma^{-1}$ by an algorithm called
standardization. This last algebra, isomorphic to Malvenuto-Reutenauer is
denoted by $\FQSym$ and is called its realization.
Similarly, the Loday-Ronco algebra of planar binary trees is isomorphic to the
algebra $\PBT$ obtained by summing the $\F_\sigma$ over sylvester classes.
Loday and Ronco then proved that $\PBT$ is isomorphic to the
free dendriform algebra on one generator.

The algebras we shall use are similar to $\FQSym$ and are indeed subalgebras
of one of its generalization, $\WQSym$. All details can be found
in~\cite{NT14}. We shall only need here their definition and simplest
properties:
define $\mFQSym$ as the algebra on elements $\F_\sigma$ indexed by
$m$-permutations whose product is given by
\begin{equation}
\F_\alpha \F_\beta = \sum_{\gamma\in\alpha\ssh\beta} \F_\gamma,
\end{equation}
where the shifted shuffle has been adapted in the following way: shift $\beta$
not by the size of $\alpha$ but only by its maximal value.
For example,
\begin{equation}
\F_{11} \F_{11} = \F_{1122} + \F_{1212} + \F_{2112} + \F_{1221} + \F_{2121}
 + \F_{2211}.
\end{equation}

It has been shown in~\cite{NT14} that there exists an analogue of $\PBT$
obtained by summing the $\F_\sigma$ again over sylvester classes, denoted by
$\mPBT$. We shall answer here the following question: is there
a natural definition of a $m$-dendriform algebra so that $\mPBT$ is
isomorphic to the free $m$-dendriform algebra on one generator?

\subsection{Dendriform algebras}

We shall first recall some highlights of the theory of dendriform algebras.
A \emph{dendriform algebra} $A$~\cite{Lod-di} is a vector space equipped with
two binary operations $\gaudend$ and $\droitdend$ respectively called left and
right, satisfying the three relations

\begin{align}
\label{un-dend}
(u\gaudend v)\gaudend w &= u\gaudend (v * w) \\
\label{de-dend}
(u\droitdend v)\gaudend w &= u\droitdend (v\gaudend w) \\
\label{tr-dend}
(u * v)\droitdend w &= u\droitdend (v\droitdend w) 
\end{align}
for any elements $u$, $v$, and $w$ in $A$
and where $*=\gaudend+\droitdend$.

It has been shown by Loday and Ronco in~\cite{LR1} that the free dendriform
algebra on one generator $D^{(1)}$ has as Poincar\'e series the generating
series of the Catalan numbers $C_n$.
We shall sketch a combinatorial proof of both this result and the fact that
$\PBT$ is isomorphic to $D^{(1)}$, since it is this proof that will be
generalized later for general values of the parameter $m$.

The proof follows three steps:
\begin{itemize}
\item The dimension of the homogeneous component of $D^{(1)}$ of degree $n$ is
at most~$C_n$.
\item Splitting the shifted shuffle product of two permutations
$\sigma\ssh\tau$ depending on its last letter being either
$\sigma_n$ or $|\sigma|+\tau_p$ endows $\FQSym$ with the structure of a
dendriform algebra.
Denote by $\DP1$ the dendriform subalgebra of $\FQSym$ generated
by $\F_1$.
\item Given any sylvester class of permutations (represented by a binary
tree), the sum of the $\F_\sigma$ over this class belongs to $\DP1$.
\end{itemize}

The first step proves that $\dim(D^{(1)}_n) \leq C_n$ whereas
the second step proves that $\dim(\DP1_n) \leq \dim(D^{(1)}_n)$ since
$\DP1$ is a dendriform algebra generated by one element, so has smaller
dimension than the \emph{free} dendriform algebra on one generator.
The third step proves that $C_n \leq \dim(\DP1_n)$ since there are
$C_n$ sylvester classes of permutations of size $n$.
All steps together prove
\begin{equation}
C_n \leq \dim(\DP1_n) \leq \dim(D^{(1)}_n) \leq C_n,
\end{equation}
hence showing both $\dim(D^{(1)}_n) = C_n$ and $\dim(\DP1_n)=\dim(D^{(1)}_n)$,
so that $\DP1$ and $D^{(1)}$ are isomorphic, $\DP1$ hence being an explicit
realization of the abstract free dendriform algebra on one generator.

\medskip
Since we shall copy these steps in the general $m$ case, let us recall
how one proves each step separately. The first step is generally proven by
orienting the dendriform relations, then showing that any expression in a
dendriform algebra is equivalent to at least a linear combination of
expressions that avoid some patterns as subexpressions, and that the total
number of expressions with $n\!-\!1$ symbols avoiding these patterns is given
by $C_n$.
It is here that we differ from the usual proof that consists in showing that
the linear combination of expressions avoiding the patterns does not depend on
the sequence of oriented relations applied to an element.
In our case, the second step consists in checking the dendriform
relations on the splitting of the associative product of $\FQSym$.
The third step amounts to see that any binary tree with respective left and
right subtrees $T_1$ and $T_2$ can be written as
$T_1\droitdend \F_1 \gaudend T_2$ if one understands $T_i$ as its expression
as an element of $\DP1$.

\section{The case $m=2$, with three operators}

Having recalled our proof in the case $m=1$, we shall detail the whole
procedure in the case $m=2$ where all ideas are already needed and where
notations and computations are easier than in the general case.

\subsection{$2$-dendriform algebras}

A $2$-dendriform algebra $A$ is a vector space equipped with three binary
operations $\gd$, $\md$, $\dd$
respectively called left, middle, and right, satisfying the
six relations

\begin{align}
\label{un-dd}
(u\gd v)\gd w &= u\gd (v * w) \\
\label{de-dd}
(u\md v)\gd w &= u\md (v\gd w) \\
\label{tr-dd}
(u\dd v)\gd w &= u\dd (v\gd w) \\
\label{qu-dd}
(u\gd v)\md w &= u\md (v\md w + v\dd w) \\
\label{ci-dd}
(u\md v+u\dd v)\md w &= u\dd (v\md w) \\
\label{si-dd}
(u *  v)\dd w &= u\dd (v\dd w) 
\end{align}
for any elements $u$, $v$, $w$ in $A$
and where $*=\gd+\md+\dd$.

\begin{note}
\label{note-rels2d}
{\rm
The operations $\gd$ and $\md+\!\dd$ satisfy the three usual
dendriform relations, so that a $2$-dendriform algebra is also a dendriform
algebra: \eqref{un-dd} is \eqref{un-dend},
adding~\eqref{de-dd} and~\eqref{tr-dd} gives back~\eqref{de-dend} while
adding~\eqref{qu-dd}, \eqref{ci-dd}, and~\eqref{si-dd} gives
back~\eqref{tr-dend}.
}
\end{note}

\begin{note}
\label{note-leroux}
{\rm
Let us compare our definition with the definition of Leroux in~\cite{Ler},
p.13, Sec. 5.1. First, it is easy to check that up to a scalar, both algebras
have a unique associative product: consider a generic element, write the
associativity axiom, and check the system involved has only one solution up to
a constant. Moreover, again in both algebras, given the operation $*$, there
are only two operations satisfying the first axiom of Leroux (which is also
Equation~\eqref{un-dd}), one of those being unsatisfactory since it is $*$
itself. So if one takes an operad morphism from the algebra of Leroux to ours,
it must, up to a scalar that can be chosen as $1$, send $*$ to $*$,
send $\gaudend$ to $\gd$, and send $\droitdend$ to $\md+\dd$.
Now, concerning $\bullet$, if one considers rule 4 of Leroux, it implies
that $\bullet$, has to be sent to our $\md$. By now, the operad morphism is
well-defined and the image of his rule 5 reads in our $3$-dendriform algebras
as
\begin{equation}
(x\md y)\md z + (x\dd y)\md z  = x\md(y\md z) + x\dd(y\md z).
\end{equation}
which is no linear combination of our relations.
}
\end{note}

\subsection{First upper bound}

Let us now define $\ddeux$ as the free $2$-dendriform algebra on one
generator. Let us also recall that the Fuss-Catalan numbers with $m=2$ are
\begin{equation}
C^{(2)}_n = \frac{1}{2n+1}\binom{3n}{n}.
\end{equation}

\subsubsection{Writing the relations as trees}

Note that any monomial in the $2$-dendriform algebra can be represented
as a complete binary tree with internal nodes labelled by one the three binary
operators and leaves by elements of the algebra.
So here follow the six $2$-dendriform relations written as trees where we have
omitted the leaves since they always are $x$, $y$, and $z$ in that order.

\begin{equation}
\label{arbs-eval-2dend}
\begin{aligned}
\treeLft{\gd}{\gd}\,&=\,\treeRgt{\gd}{*}
&
\treeLft{\gd}{\md}\ &=\ \treeRgt{\md}{\gd}
&
\treeLft{\gd}{\dd}\ &=\ \treeRgt{\dd}{\gd}
\\
\treeLft{\md}{\gd}\ &=\ \treeRgt{\md}{\md} \!\!\!+\!\!\! \treeRgt{\md}{\dd}
\quad &\quad
\treeLft{\md}{\md} \!\!\!+\!\!\! \treeLft{\md}{\dd}\ &=\ \treeRgt{\dd}{\md}
\quad &\quad
\treeLft{\dd}{*}\ &=\ \treeRgt{\dd}{\dd}
\end{aligned}
\end{equation}

Note that the trees with the symbol $*$ represent the sum of all trees where
$*$ is replaced by the three symbols $\gd$, $\md$, and $\dd$.

Let us now prove

\begin{theorem}
\label{th-ddeux}
Let $n$ be a positive integer. Then
\begin{equation}
\dim(\ddeux_n) \leq C^{(2)}_n.
\end{equation}
\end{theorem}

\subsubsection{Orienting the relations}

Let us orient the $2$-dendriform relations as follows:
\begin{equation}
\begin{aligned}
\treeLft{\gd}{\gd}\,&\rightarrow \ \treeRgt{\gd}{*}
&
\treeLft{\gd}{\md}\ &\leftarrow \ \treeRgt{\md}{\gd}
&
\treeLft{\gd}{\dd}\ &\leftarrow \ \treeRgt{\dd}{\gd}
\\
\treeLft{\md}{\gd}\ &\rightarrow \
                     \treeRgt{\md}{\md} \!\!\!+\!\!\! \treeRgt{\md}{\dd}
\quad &\quad
\treeLft{\md}{\md} \!\!\!+\!\!\! \treeLft{\md}{\dd}\
                    &\leftarrow \ \treeRgt{\dd}{\md}
\quad &\quad
\treeLft{\dd}{*}\   &\leftarrow \ \treeRgt{\dd}{\dd}
\end{aligned}
\end{equation}

This amounts to forbid the following six tree patterns:
\begin{equation}
\label{6arbs}
\treeLft{\gd}{\gd}
\qquad
\treeRgt{\md}{\gd}
\qquad
\treeRgt{\dd}{\gd}
\qquad
\treeLft{\md}{\gd}
\qquad
\treeRgt{\dd}{\md}
\qquad
\treeRgt{\dd}{\dd}
\end{equation}

To get to the conclusion that $\dim(\ddeux)_n \leq C^{(2)}_n$, we need to
prove two facts: first, orienting the $2$-dendriform relations as before
brings ultimately linear combinations of trees avoiding as subtrees the six
trees presented in~\eqref{6arbs}; and the number of binary trees
avoiding these six tree patterns is enumerated by $C^{(2)}_n$.
We say here that a tree $T$ avoids a subtree of size $2$ if there is no edge
oriented in the same way and with the same parent and child as in the subtree.

\subsubsection{Rewriting trees with forbidden patterns}

\begin{proposition}
Any tree can be expressed as a linear combination of trees avoiding the six
forbidden tree patterns.
\end{proposition}

\Proof
Let us consider a tree with $n$ operators.
If $n\leq2$, the statement holds.

Otherwise, we shall prove the statement by induction. Assume that any tree
with at most $n-1$ operators satisfies the statement.
Hence, given a tree $T$ with subtrees $T_1$ and $T_2$, one can assume
that both $T_1$ at $T_2$ avoid the six forbidden trees since one can rewrite
separately the left and right subtrees of $T$.
Let us now split the cases depending on the root of $T$.

First, if the root of $T$ is $\gd$, then $T$ avoids the six patterns if the
root of $T_1$ is not $\gd$, so we can assume that the root of $T_1$ is $\gd$
and denote by $T'_1$ and $T''_1$ its own subtrees. We can then apply the
rewriting

\begin{equation}
T =
\treeGen{\gd}{T_1}{T_2} =
\treeGenG{\gd}{\gd}{T'_1}{T''_1}{T_2} =
\treeGenD{\gd}{T'_1}{*}{T''_1}{T_2}.
\end{equation}
Note that since $T_1$ avoids all patterns, the root of $T'_1$ is different
from $\gd$, so the patterns can only be found in the right subtree with root
label $*$. But, again by induction, this tree is a linear combination of trees
avoiding the patterns, and so is the linear combination obtained for $T$.

Now, if the root of $T$ is $\md$, then either the root of $T_1$ is $\gd$ or
is not. If the root of $T_1$ is not $\gd$, then $T$ avoids the six patterns if
the root of $T_2$ is not $\gd$, so we can assume that the root of $T_2$ is
$\gd$ and denote by $T'_2$ and $T''_2$ its subtrees.
We can then apply the rewriting 

\begin{equation}
T =
\treeGen{\md}{T_1}{T_2} =
\treeGenD{\md}{T_1}{\gd}{T'_2}{T''_2} =
\treeGenG{\gd}{\md}{T_1}{T'_2}{T''_2}.
\end{equation}
Since $T_2$ avoids all patterns, the root of $T'_2$ is not $\gd$, and
this last tree avoids all patterns.

Now, assume that the root of $T_1$ is $\gd$. In that case, we
have the following picture:
\begin{equation}
T =
\treeGen{\md}{T_1}{T_2} =
\treeGenG{\md}{\gd}{T'_1}{T''_1}{T_2} =
\treeGenD{\md}{T'_1}{\md}{T''_1}{T_2} +
\treeGenD{\md}{T'_1}{\dd}{T''_1}{T_2} .
\end{equation}
Since $T_1$ avoids all patterns, the root of $T'_1$ is different
from $\gd$, and the two trees on the right belong to a case seen previously
where the root is $\md$ and the root of its left subtree is not $\gd$. So by
induction, all these trees rewrite as linear combinations of trees avoiding
the six patterns.

The case where the root is a $\dd$ is analogous to the previous cases and is
treated without difficulties.
\qed

\subsubsection{Enumerating the trees avoiding the patterns}

The technique to enumerate a set of trees with forbidden patterns is always
the same. Define by $U$ the generating series of all these trees enumerated
thanks to their number of operators, and by $U_\gd$, $U_\md$, and $U_\dd$ the
subseries of these trees where the root respectively is $\gd$, $\md$, and
$\dd$.
We then have the following system of equations directly derived from the
forbidden patterns:
\begin{equation}
\label{eqs-forb}
\left\{
\begin{aligned}
U     &= 1 + U_\gd + U_\md + U_\dd \\
U_\gd &= x\, (U-U_\gd)\, U \\
U_\md &= x\, (U-U_\gd)^2 \\
U_\dd &= x\, U
\end{aligned}
\right.
\end{equation}

One then easily eliminates the three partial series and gets
\begin{equation}
U = (1+ x U)^3,
\end{equation}
which is the generating series of non-empty ternary trees enumerated by their
number of nodes minus $1$.
Since $\dim(D^{(2)}_n)$ is smaller than or equal to the coefficient of
$x^{n-1}$ in $U$ (all operators are binary so there is one more operand than
operators), we conclude to the theorem:
$\dim(D^{(2)}_n) \leq C^{(2)}_{n}$.

\begin{note}
{\rm
One needs to carefully choose the forbidden patterns since, choosing
\emph{e.g.}, to forbid the other pattern of~\eqref{de-dd}, one gets a
non-optimal hence useless upper bound.
}
\end{note}

\subsection{An explicit $2$-dendriform algebra}

This research has of course been done in a different order than the
presentation given here. To get the correct rewriting rules with no explicit
object would have been pointless and it has indeed been done the other way
round: first guess a way to split the shuffle product of $2$-permutations into
three operations that separate their sylvester classes. And then look
for the relations these operations satisfy. We shall now present how we
adressed this question.

When one computes $\F_{11}^2$, one finds three sylvester classes, and more
precisely, splitting according to the dendriform operations $\gaudend$ and
$\droitdend$, one has

\begin{align}
\label{gaud11}
\F_{11} \gaudend \F_{11}   &= \F_{1221} + \F_{2121} + \F_{2211}, \\
\label{rigd11}
\F_{11} \droitdend \F_{11} &= \F_{1122} + \F_{1212} + \F_{2112}.
\end{align}

The first set of $2$-permutations constitutes a single sylvester class,
whereas the second one can be split in two, hence suggesting that one must
split $\droitdend$ into two operations, and preserve $\gaudend$. Hence
justifying Note~\ref{note-rels2d}.

\subsubsection{Splitting $*$ into three operations}

Let us now get to the construction itself.
We shall use the right-to-left minima of $2$-permutations
$\tau=\tau_1\dots\tau_p$
that is the set of values $i\in [1,p]$ such that
$\tau_j\geq\tau_i$ for all $j>i$.
For example, the right-to-left minima of $212313$ are in decreasing order
$\{6,5,2\}$ whereas the right-to-left minima of $4121235453$ are
$\{10,6,5,4,2\}$.
We shall write $m_1(\tau)$ for the first right-to-left minimum, that is always
$p$, and $m_2(\tau)$ for the second one.

Write $\tau$ as $\tau = \tau' \cdot \tau'' \cdot \tau_p$ where the first
letter of $\tau''$ is $m_2(\tau)$.
Note that any $2$-permutation has at least \emph{two} right-to-left minima
thanks to the definition, so this way of splitting $\tau$ always makes
sense.

For example,
\begin{equation}
212313 = 2123.1.3
\qquad\text{and}\qquad
4121235453 = 41212.3545.3.
\end{equation}

Then define the three operations on
$\sigma=\sigma_1\dots\sigma_n$
and $\tau[n]=\tau'\cdot \tau'' \cdot \tau_p$ as
\begin{equation}
\left\{
\begin{aligned}
\sigma\gd\tau &= (\sigma_1\dots\sigma_{n-1} \shuffle \tau[n])
                 . \sigma_n \\
\sigma\md\tau &= (\sigma\shuffle (\tau'\cdot\tau'')) . \tau_p
                 - \sigma\dd\tau \\
\sigma\dd\tau &= (\sigma\shuffle\tau') . \tau''.\tau_p
\end{aligned}
\right.
\end{equation}

In other words, $\sigma\gd\tau$ (respectively $\sigma\md\tau$ and
$\sigma\dd\tau$) is the subset of $\sigma\ssh\tau$ where the right-most letter
of $\sigma$ ends up after $m_1(\tau)$ (resp. between $m_2(\tau)$ and
$m_1(\tau)$, resp. before $m_2(\tau)$).

For example,
\begin{equation}
11\md 11 = 1212 + 2112
\qquad\text{and}\qquad
11\dd 11 = 1122 .
\end{equation}

So this splits the right hand-side of~\eqref{rigd11} into its two
sylvester classes.
As a bigger example, one also has
\begin{equation}
\begin{split}
11\md 2112 &= 132213 + 312213 + 321213 + 322113, \\
11\dd 2112 &= (11\shuffle 32) . 23 \\
           &= 113223+131223+132123+311223+312123+321123.
\end{split}
\end{equation}

\subsection{The $2$-dendriform structure}

Let us now check that the three operations defined in the last paragraph
indeed define a $2$-dendriform algebra, that is, satisfy
Relations~\eqref{un-dd} up to~\eqref{si-dd}.

Relation~\eqref{un-dd} is automatic since it is a known dendriform
relation. Relations~\eqref{de-dd} and~\eqref{tr-dd} sum up to another known
dendriform relation, so we only have to check one relation to get the other.
Relation~\eqref{tr-dd} is indeed true since the values of the right-to-left
minima of $v$ are equal to the values of the right-to-left minima of $v\gd w$:
in particular, $m_1(v)=m_1(v\gd w)$ and $m_2(v)=m_2(v\gd w)$.
So $u$ is shuffled with the same prefix of $v$ on both sides of the relation
regardless of $w$.

The last three relations all sum up to the last dendriform relation, so we
only need to prove two of those to get the last one. Relation~\eqref{si-dd} is
direct since it amounts in both terms to compute $(u\ssh v\ssh w').w''$ where
$w=w'.w''$ is the splitting of $w$ such that $w''_1=m_2(\tau)$.
Relation~\eqref{qu-dd} is also true since in both expressions the sum consists
in the elements such that the last letter of $u$ is to the right
of $v$ and between the first two right-to-left minima of $w$. Indeed, the
second right-to-left minimum of $v\md w + v\dd w$ is always to the right of
the last letter of $v$ and hence is the second right-to-left minimum of $w$.

We can then conclude:

\begin{theorem}
\label{fqs-2d}
The three operations defined previously endow ${}^2\FQSym$ with the structure
of a $2$-dendriform algebra.
\end{theorem}

We shall denote by $\DP2$ the $2$-dendriform subalgebra of ${}^2\FQSym$
generated by the element $\F_{11}$.

\begin{note}
\label{not-2fq-2dend}
{\rm
The previous theorem shows in particular that the dimension of the homogeneous
component of size $n$ of $\DP2$ is smaller than or equal to the dimension of
the same component of $\D^{(2)}$.
}
\end{note}

\subsection{$2$-dendriform products and sylvester classes}
\label{sylv-2d}

Let us now prove that the three $2$-dendriform products satisfy the property
they were designed for, that is, send a pair of sylvester classes to an union
of sylvester classes.
First, if the product is $\gd$, it comes from known results on the dendriform
structure of ${}^2\FQSym$, see~\cite{NT14}. The same argument holds if one
considers $\md+\dd$ so we only need to prove that a product, say $\dd$, of
two sylvester classes is indeed an union of sylvester class.

Let us then consider two words $w$ and $w'$ that are equivalent up to one
sylvester rewriting: up to exchanging the roles of $w$ and $w'$, we have
$w=\dots ac\dots b\dots$ and $w'=\dots ca\dots b\dots$ with $a\leq b<c$.
Assume now that $w$ belongs to a product $u\dd v$.

If the three letters $a$, $b$, $c$ belong to $u$, then $w'$ 
belongs to $u'\dd v$ where $u$ and $u'$ are in the same sylvester class.
The same holds if the three letters belong to $v$. Now, if the three letters
do not belong to the same word, then $a$ belongs to $u$ and $c$ belongs to
$v$, and $w'$ also belongs to $u\dd v$: 
if $b$ belongs to $u$, $a$ is not the right-most letter of $u$ and hence is
swapped with $c$ in $u\dd v$ and if $b$ belongs to $v$, then $c$ cannot be
$m_2(v)$ since $c>b$ so $a$ can again be swapped with $c$ in $u\dd v$.

So since $\DP2$ is generated by $\F_{11}$ which is a sylvester class all by
itself, $\DP2$ has a basis given by linear combinations of sylvester classes
of $2$-permutations which is coherent with the combination of
Note~\ref{not-2fq-2dend} and Theorem~\ref{th-ddeux} that together show that
$\dim\DP2_n\leq C_n^{(2)}$.
We shall see in the next paragraph that any sylvester class belongs to $\DP2$.

\subsection{Sylvester classes of $2$-permutations as $2$-dendriform products}

Let us summarize what we have proved up to now.
We showed in Theorem~\ref{th-ddeux} that the free $2$-dendriform algebra on
one generator has dimension at most $C^{(2)}_n$ in size $n$.
We showed in Theorem~\ref{fqs-2d} that our way of cutting the shifted shuffle
product of permutations in three has indeed the structure of a $2$-dendriform
algebra, so that (Note~\ref{not-2fq-2dend}) the dimension of the
$2$-dendriform algebra $\DP2$ generated by $\F_{11}$ is at most $C^{(2)}_n$ in
size $n$ too.
Finally, we showed in Section~\ref{sylv-2d} that $\DP2$ has a basis given by
linear combinations of sylvester classes of $2$-permutations.

Since sylvester classes of $2$-permutations are enumerated by the Fuss-Catalan
numbers, we only have to prove that any sylvester class of $2$-permutations
belong to $\DP2$ to conclude. Indeed, it would prove that $\dim(\DP2)_n$
is of dimension at least $C^{(2)}_n$, hence of dimension $C^{(2)}_n$ exactly,
so that $\dim(D^{(2)}_n)$ also is the same.

This last result is indeed correct and comes from the following three
relations we shall prove in the next section.
The relations are written in terms of ternary trees since sylvester classes of
$2$-permutations are indeed encoded by ternary trees (see
Section~\ref{sec-sylvmp}).
\begin{equation}
\label{gdgd}
\treeGenT{\bulle}{T_1}{T_2}{T_3} =
\treeGenTd{\bulle}{T_1}{T_2} \gd T_3.
\end{equation}
\begin{equation}
\label{mdmd}
\treeGenTs{\bulle}{T'_1}{\bulle}{T_1}{T_2}{T'_2} =
\treeGenTd{\bulle}{T_1}{T_2} \md
\treeGenTd{\bulle}{T'_1}{T'_2}.
\end{equation}
\begin{equation}
\label{dddd}
\treeGenTz{\bulle}{T} = T \dd \bulle\,.
\end{equation}

\subsection{Proving the three reduction relations}

Let us first show that all these equations are correct. We shall see that
the products on the right hand-side always consist in exactly one tree and
that it is the tree given on the left hand-side.

Indeed, let us translate the trees as their canonical words obtained by
reading the unique labelling of their tree as ternary search tree (see
Section~\ref{sec-sylvmp} for details about this construction).
Then, since canonical words of sylvester classes are the words avoiding the
pattern $132$, a product of two sylvester classes, being an union of sylvester
classes, contains as many terms as there are words avoiding $132$ in the
shifted shuffle of their canonical elements.
So one can restrict to the corresponding product of canonical words.

For example, the product in~\eqref{gdgd} reads as $T_1 rT_2r \gd T_3$ where
$r$ is the value of the root. The product is $(T_1 rT_2\shuffle T_3[n]) r$ if
$n$ denotes the number of nodes of the tree on the left.
Since all letters of $T_3[n]$ are greater than both $T_1$ and $r$, and that
all letters of $T_1$ are smaller than $r$, there is only one word avoiding
$132$ in this set. It is $T_3[n] T_1 r T_2 r$, which is exactly the
canonical word of the tree on the left hand-side of the equation.

The product in~\eqref{mdmd} is proven in the same way: one has to compute
$T_1 rT_2r \md T'_1r'T'_2r'$. This product is equal to
$(T_1 rT_2 \shuffle (T'_1r'T'_2)[n]) r r'[n]$. Again all letters of
$(T'_1r'T'_2)[n]$ are greater than $r$ whereas all letters of $T_1 rT_2$ are
smaller than or equal to $r$, so that there is only one word avoiding $132$ in
this shuffle product, which is
$T'_1 r' T'_2 T_1 r T_2 r r'$.
Equation~\eqref{dddd} is immediate.

\subsection{Concluding remarks}

Thanks to the three relations given before, one now sees that any
sylvester class is obtained as a suitable product of the generator $\F_{11}$
by induction on the sizes of the trees.
So the $2$-dendriform algebra generated by $\F_{11}$ has a basis indexed
by sylvester classes of $2$-permutations.

For example,
\begin{equation}
\begin{split}
\treeEx &= \treeExb \gd \bulle
         = \left(\treeExc \md \treeExd\right) \gd \bulle \\
        &= \left( \left(\bulle \md \treeExe\right) \md \treeExd
           \right) \gd \bulle \\
        &= \left( \left(\bulle \md \left( \treeExf \dd \bulle\right) \right)
           \md \treeExd \right) \gd \bulle \\
        &= \left( \left(\bulle \md \left( \left(\bulle\md\treeExg\right)
           \dd \bulle\right) \right)
           \md \left(\bulle \md \treeExh \right) \right) \gd \bulle \\
        &= \left( \left(\bulle \md \left(
                  \left(\bulle\md (\bulle\md\bulle) \right)
           \dd \bulle\right) \right)
           \md \left(\bulle \md (\bulle\dd\bulle) \right) \right)
           \gd \bulle 
\end{split}
\end{equation}

We can summarize all the previous results as

\begin{theorem}
The free $2$-dendriform algebra $D^{(2)}$ has as series of dimensions the
series of Fuss-Catalan numbers $C^{(2)}_n$ and the $2$-dendriform algebra
$\DP2$ that is a $2$-dendriform subalgebra of ${}^2\FQSym$ is free and
isomorphic to $D^{(2)}$.
\end{theorem}

\section{The general case}

We shall now consider the case where $m$ is any integer and not only $1$ or
$2$. Those two cases were discussed before and the general case is just an
adaptation of those for general $m$.

The strategy is once again the same:
\begin{itemize}
\item define an abstract $m$-dendriform algebra by providing relations
between $m\!+\!1$ operations,
\item prove that the Fuss-Catalan numbers provide an upper bound of the
dimension of the homogeneous component of size $n$ of the free $m$-dendriform
algebra on one generator $D^{(m)}$,
\item define a way to split the shifted shuffle of $m$-permutations into
$m\!+\!1$ operations and check that these operations satisfy the
$m$-dendriform relations,
\item show that these $m\!+\!1$ operations send pairs of sylvester classes
into unions of sylvester classes,
\item show conversely that each sylvester class can be obtained as a suitable
product from the generator $\F_{1^m}$.
\end{itemize}

Let us now show how each step presents itself in the generic case.

\subsection{$m$-dendriform algebras}

An \emph{$m$-dendriform algebra} is an algebra with $m+1$ operations labelled
$\gd$, $\md_1,\dots,\md_{m-1}$, $\dd$ subject to the relations

\begin{align}
\label{un-md}
(u\gd v)\gd w &= u\gd (v * w), \\
\label{de-md}
(u \md_i v)\gd w &= u \md_i (v\gd w),
& \text{for all $i$,} \\
\label{tr-md}
(u \dd v)\gd w &= u \dd (v\gd w), \\
\label{qu-md}
(u\gd v) \md_i w &= u \md_i (v\dd w + \sum_{j\geq i} v \md_j w),
& \text{for all $i$,} \\
\label{ci-md}
(u * v) \dd w &= u \dd (v \dd w), \\
\label{si-md}
(u\dd v + \sum_{j\geq m-i} u \md_j v) \md_i w &= u \dd (v \md_i w),
& \text{for all $i$}, \\
\label{se-md}
(u \md_k v) \md_i w &= u \md_{k+i} (v \md_i w),
& \text{if $k+i<m$,}
\end{align}
for all $u$, $v$, $w$, and
where $*=\gd+\md_1+\dots+\md_{m-1}+\gd$.

Note that all expressions with two rules are represented in the products so
that adding all these relations together shows that $*$ is associative.

\begin{note}
{\rm
There are respectively $1$, $m-1$, $1$, $m-1$, $1$, $m-1$ and
$\binom{m-1}{2}$ relations of the corresponding types hence $\binom{m+2}{2}$
relations in total, so that, starting with $2(m+1)^2$ expressions involving
two operators, there only remains $2(m+1)^2-\binom{m+2}{2}$ different
expressions,
which is
the pentagonal number $(m+1)(3m+2)/2=\binom{3m+3}{2}/3$,
which is exactly the third Fuss-Catalan number
$C^{(m)}_3=\frac{1}{3m+1}\binom{3m+3}{3}$.
}
\end{note}

\begin{note}
{\rm
Thanks to this definition, any $m$-dendriform algebra is also an
$(m\!-\!1)$-dendriform algebra with the operations
$\gd$, $\md_1,\dots,\md_{m-2}$, and $\md_{m-1}+\dd$.
Indeed, the relations concerning this last operation are obtained by summing
together relations involving both $\md_{m-1}$ and $\dd$.
}
\end{note}

\begin{note}
\label{note-leroux-gen}
{\rm
As it was already the case with $m=2$, our algebras are not isomorphic to the
algebras of Leroux defined in~\cite{Ler}. Indeed, the same technique as in
Note~\ref{note-leroux} proves the fact.
}
\end{note}

\subsection{The upper bound}

Once again, we orient the relations so that, the forbidden tree patterns are
\begin{equation}
\treeLft{\gd}{\gd}
\qquad
\treeRgt{\md_i}{\gd}
\qquad
\treeRgt{\dd}{\gd}
\qquad
\treeLft{\md_i}{\gd}
\qquad
\treeRgt{\dd}{\dd}
\qquad
\treeRgt{\dd}{\md_i}
\end{equation}
for all $i$
and all trees
\begin{equation}
\treeRgt{\md_j}{\md_i}
\end{equation}
with $1 \leq i < j \leq m-1$.

\subsubsection{Removing the forbidden tree patterns}
\label{sub-forbid}

\begin{proposition}
Any tree can be expressed as a linear combination of trees avoiding the
general forbidden tree patterns.
\end{proposition}

The proof is the exact copy of the proof in the case $m=2$: prove it by
induction on trees, and split the cases according to the operator on the
root. Then solve the cases starting with $\gd$, then all $\md_i$ in increasing
order then $\dd$. This order guarantees that given a tree that has a forbidden
pattern, either its rewriting does not have any pattern anymore, or the
pattern has already been dealt with before. Indeed, this works because the
forbidden patterns were not chosen at random. They all are the patterns in
their corresponding equation that are maximal in a certain sense:
they are the patterns whose root is maximal according to the natural order of
the operators and whose other vertex is extremal (maximal if the root is $\dd$
and minimal otherwise). So their rewritings are smaller in a certain sense so
that the algorithm necessarily ends.

\subsubsection{Enumerating the trees avoiding the patterns}

Let us again denote by $U$ the generating series of all the trees avoiding the
forbidden patterns enumerated thanks to their number of operators, and by
$U_\gd$, $U_{\md_i}$, $U_\dd$ the subseries of these trees where the root
respectively is $\gd$, $\md_i$, or $\dd$.
We then have the following system of equations directly derived from the
forbidden patterns:
\begin{equation}
\label{eqs-forb-gen}
\left\{
\begin{aligned}
U     &= 1 + U_\gd + \sum_{i} U_{\md_i} + U_\dd\\
U_\gd &= x\, (U-U_\gd)\, U \\
U_{\md_1} &= x\, (U-U_\gd)^2 \\
U_{\md_2} &= x\, (U-U_\gd)(U-U_\gd-\sum_{i<2} U_{\md_i}) \\
\vdots\quad  &= \qquad\qquad \vdots \\
U_{\md_{m-1}} &= x\, (U-U_\gd)(U-U_\gd-\sum_{i<m-1} U_{\md_i}) \\
U_{\dd} &= x\, U . 1 \\
\end{aligned}
\right.
\end{equation}
Then one easily checks that this rewrites as
$U_\gd = U \frac{xU}{1+xU}$,
and
$U_{\md_i} = U \frac{xU}{(1+xU)^{i+1}}$
for all $i$, and $U_{\dd}=xU$,
so that one directly obtains
\begin{equation}
\begin{split}
U &= 1+xU + U\left(\sum_{i=0}^{m-1}\frac{xU}{(1+xU)^{i-1}}\right) \\
  &= 1+xU + U\left(1-\frac{1}{(1+xU)^{m}}\right)
\end{split}
\end{equation}
which implies
\begin{equation}
U = (1+xU)^{m+1},
\end{equation}
hence the same conclusion as in the case $m=2$ since this series enumerates
the $(m+1)$-ary trees by their number of nodes minus $1$.

Note that the general computation of $U$ gives away how to choose the
forbidden patterns.

\subsection{An explicit $m$-dendriform algebra: Splitting the shifted shuffle
into $m+1$ parts}

The recipe is the same as in the $m=2$ case:
given an $m$-permutation, compute its first $m$ right-to-left minima.
Then define $u \md_i v$ as the subset of $u\ssh v$ where the last letter of
$u$ is between $m_{i+1}(v)$ and $m_i(v)$.
In particular, $\md_1$ is equal to the operator $\md$ of the $2$-dendriform
algebra and $\dd$ just selects the words where the last letter of $u$ is to
the left of $m_m(v)$.

One then easily checks that this way of splitting the shifted shuffle indeed
endows $\mFQSym$ of a $m$-dendriform algebra structure since all
relations hold.
Relation~\eqref{un-md} holds since it is a dendriform relation.
Relations~\eqref{de-md} and~\eqref{tr-md} hold for the same reason
Relations~\eqref{de-dd} and~\eqref{tr-dd} hold in the $2$-dendriform algebra:
the left-to-right minima of $v\gd w$ are the same as the left-to-right minima
of $v$.
Relation~\eqref{qu-md} holds since the two conditions: the last letter of $u$
is to the right of the last letter of $v$ and the last letter of $u$ is
between $m_{i+1}(w)$ and $m_i(w)$ are equivalent to the two conditions:
the last letter of $v$ is to the left of $m_i(w)$ and
the last letter of $u$ is between the same letters $m_{i+1}(w)$ and $m_i(w)$
and to the right of the last letter of $v$.

Relation~\eqref{ci-md} holds for the same reason as in the $2$-dendriform
algebra.
Relation~\eqref{si-md} holds since the two conditions: the last letter of $u$
is to the left of $m_{m-i}(v)$ and the last letter of $v$ is to the left of
$m_{i}(z)$ are equivalent to the two conditions:
the last letter of $u$ is to the left of $m_{m}(v\md_i w)$ and the last
letter of $v$ is to the left of $m_{i}(z)$.
Finally, Relation~\eqref{se-md} holds again by equivalence of the two
pairs of conditions.

\subsection{Multiplying sylvester classes}

Proving that any product of sylvester classes is an union of sylvester classes
is done exactly as in the $m=2$ case: either both letters $a$ and $c$ come
from the same word and it is obvious, either they do not but in that case,
either $b$ was with $a$ and one can exchange $a$ and $c$, or $b$ was with $c$
but then $c$ cannot be a right-to-left minimum, hence can also be exchanged
with $a$.

\subsection{Obtaining a particular sylvester class}

As in the $m=2$ case, there exists $m+1$ relations that allow one to
strip a tree off its $m+1$ subtrees, in order. The first and last step of this
prodedure are as easy as in the $m=2$ case so we shall concentrate on the case
where we want to obtain a tree whose $k$ subtrees from the right are empty.
It just amounts to the following formula:
\begin{equation}
\treeGenTes{\bulle}{T'_1 \dots\ \ \ \ \ }{T'_{m-k}\ }%
{\bulle}{T_1\dots\ \ }{T_m\ \ }{T'_{m+1-k}} =
\treeGenTd{\bulle}{T_1\dots}{T_m}
\md_k
\treeGenTd{\bulle}{T'_1\dots\ \ \ \ \ \ }{T'_{m+1-k}}.
\end{equation}
  
This formula holds since it amounts to compute
$T_1r \dots r T_m r \md_k T'_1r'\dots T'_{m+1-k}r'^{k}$
which is equal to
\begin{equation}
(T_1r \dots r T_m \shuffle (T'_1r'\dots T'_{m+1-k})[n]) r r'^k,
\end{equation}
where the left (respectively right) part of the shuffle consists in letters
smaller (resp. greater) than $r$, so that there is only one word in this
shuffle product that avoids the pattern 132 and it is
\begin{equation}
(T'_1r'\dots T'_{m+1-k})[n])\ T_1r \dots r T_m\  r r'^k,
\end{equation}
which is the corresponding $132$-avoiding $m$-permutation of the tree on the
left of the formula.

Piecing all the results of this Section together, we get
\begin{theorem}
The free $m$-dendriform algebra $D^{(m)}$ has as series of dimensions the
series of Fuss-Catalan numbers $C^{(m)}_n$ and the $m$-dendriform subalgebra
$\DP m$ of $\mFQSym$ generated by $\F_{1^m}$ is free and isomorphic to
$D^{(m)}$.
\end{theorem}

\section{The dual of the $m$-dendriform algebras}

The dual of the dendriform algebras is well-known and is the $Dias$ algebras,
some monoidal algebras.
Since our relations are very similar to the dendriform relations, it shall not
come as a surprise that the $m$-dendriform algebras also have duals and that
these duals generalize $Dias$.

Recall that the dual of a quadratic algebra $A$ is defined as an algebra on
the same number of generators whose relations are orthogonal to the relations
of $A$. Here, since the relations defining the $m$-dendriform algebras do not
share any $2$-vertices tree, the computation of the orthogonal is as easy as
in the dendriform case.

\subsection{Presentation of $m$-Dias}

The dual of the $m$-dendriform algebras defined as $m$-Dias algebras are
defined as algebras on $m\!+\!1$ binary operations labelled $\dgd$, $\dmd_1$,
$\dots$, $\dmd_{m-1}$, $\ddd$ satisfying the set of relations

\begin{align}
\label{un-dmd}
(u\dgd v)\dgd w &= u\dgd (v\dgd w) = u\dgd (v\dmd_i w) = u\dgd (v\ddd w),
& \text{for all $i$,} \\
\label{de-dmd}
(u \dmd_i v)\dgd w &= u \dmd_i (v\dgd w),
& \text{for all $i$,} \\
\label{tr-dmd}
(u \ddd v)\dgd w &= u \ddd (v\dgd w), \\
\label{qu-dmd}
(u\dgd v) \dmd_i w &= u \dmd_i (v\ddd w) = u \dmd_i (v \dmd_j w),
& \text{for all $j\geq i$,} \\
%
\label{ci-dmd}
(u \dgd v) \ddd w &= (u \dmd_i v) \ddd w = (u \ddd v) \ddd w
                  = u \ddd (v \ddd w),
& \text{for all $i$,} \\
\label{si-dmd}
(u\dmd_j v)\dmd_i w &= (u\ddd v)\dmd_i w =  u \ddd (v \dmd_i w),
& \text{for all $i\!+\!j\!\geq\! m$}, \\
\label{se-dmd}
(u \dmd_k v) \dmd_i w &= u \dmd_{k+i} (v \dmd_i w),
& \text{for all $k\!+\!i<m$,}
\end{align}
for all words $u$, $v$, and $w$.

\subsection{The free $m$-Dias algebra on one generator}

As is customary for the dual of quadratic algebras, one obtains the elements
of the dual by considering the trees formed only by the forbidden tree
patterns of the original algebra. It amounts here to consider the patterns
\begin{equation}
\treeLft{\dgd}{\dgd}
\qquad
\treeRgt{\dmd_i}{\dgd}
\qquad
\treeRgt{\ddd}{\dgd}
\qquad
\treeLft{\dmd_i}{\dgd}
\qquad
\treeRgt{\ddd}{\ddd}
\qquad
\treeRgt{\ddd}{\dmd_i}
\end{equation}
for all $i$ and all trees
\begin{equation}
\treeRgt{\dmd_j}{\dmd_i}
\end{equation}
with $1 \leq i < j \leq m-1$.
We shall call these patterns the \emph{valid} patterns in the sequel.

As before, we shall prove that the free $m$-Dias algebra on one
generator, denoted here by $Di^{(m)}$ satisfies
\begin{equation}
\dim Di^{(m)}_n = \binomial{n+m-1}{n-1}.
\end{equation}

Note that this number is equal to the number of trees having only valid tree
patterns. Indeed, denoting again by $U$ and $U_{\bullet}$ the total number of
trees or the number of trees with a given root, we have
\begin{equation}
\label{eqs-ok-gen}
\left\{
\begin{aligned}
U     &= 1 + U_\dgd + \sum_{i} U_{\dmd_i} + U_\ddd\\
U_\dgd &= \frac{x}{1-x} \\
U_{\dmd_1} &= x\, (1+U_\dgd)^2 \\
U_{\dmd_2} &= x\, (1+U_\dgd)(1+U_\dgd+\sum_{i<2} U_{\dmd_i}) \\
\vdots\quad  &= \qquad\qquad \vdots \\
U_{\dmd_{m-1}} &= x\, (1+U_\dgd)(1+U_\dgd+\sum_{i<m-1} U_{\dmd_i}) \\
U_{\ddd} &= x\, U \\
\end{aligned}
\right.
\end{equation}
Then one easily checks that this rewrites as
$U_{\dmd_i} = \frac{x}{(1-x)^{i+1}}$
for all $i$, and $U_{\dd}=xU$,
so that one directly obtains
\begin{equation}
U = 1+xU + \sum_{i=1}^{m}\frac{x}{(1-x)^{i}}
\end{equation}
so that $U=\frac{1}{(1-x)^{m+1}}$.

Let us observe that the $m$-Dias algebras are monoidal (all relations involve
one tree equal to another). So we have to prove that any tree is equivalent up
to the relations to exactly one tree that has only the valid tree patterns.

\subsubsection{Upper bound}

First, let us see why any tree is equivalent to a tree with only valid tree
patterns. The method is the same as for $m$-dendriform algebras: do it by
induction on trees, splitting the cases according to the operator on the
root.
Indeed, the same argument as presented in the $m$-dendriform case (see
Section~\ref{sub-forbid}) works since the valid patterns, being the same as
the forbidden patterns in the $m$-dendriform case are \emph{maximal} in a
certain sense among their rewritings. So all trees rewrite to smaller trees in
a certain sense, hence get equal to a tree with only valid patterns by
induction on that particular order.

\subsubsection{Lower bound}

We now want to show that two trees with only valid patterns cannot be
rewritten into one another, that is, do not belong to the same class of trees.
Let us explain how this proof works using rewriting graphs and a confluence
property. So let us consider the graph $G$ of all rewritings between $2$-edges
subtrees and their corresponding non-forbidden $2$-edges subtree in the
$m$-Dias algebras. The trees in the connected components of $G$ are all trees
that are equivalent in the sense of the previous section. 
Now consider $G'$ the graph obtained from $G$ by orienting all edges of $G$
by rewriting towards the valid subtrees.
Note that in our case $G'$ is acyclic since our valid subtrees, as already
noted in Section~\ref{sub-forbid} are maximal in a certain sense. We shall
assume here that $G'$ is acyclic.

\medskip
Let us now assume that there is a connected component of $G$ containing
strictly more that one tree with only valid patterns.
Then inside $G'$, consider the subcomponents going to each of the trees with
only valid patterns. These subcomponents cannot be nonintersecting: since $G$
is connected, there has to be two trees $T$ and $T'$ that are related by an
edge of $G$ so that $T$ and $T'$ belong to different subcomponents. Then,
depending on the orientation of the rewriting, either $T$ or $T'$ belongs to
both subcomponents. So, there exists a tree that can go in $G'$ to two
different trees with only valid patterns.
Then consider such a tree $T$ minimal in the sense that it has two successors
going towards different valid trees that do not share the previous property.
Such a tree $T$ exists since $G'$ is acyclic.
Then these two successors have disjoint sets of images in $G'$.

So we shall here show that given a tree $T$ and any two edges in $G'$ starting
from $T$ defining two trees $T'$ and $T''$, both trees have images in common
in $G'$ hence proving that all connected components of $G$ have only one tree
with only valid patterns by this \emph{confluence} property. Note that if the
two rewritings do not share a vertex, the property is obvious. So we only have
to prove the property on rewritings having a common vertex, hence on trees
with three vertices. Since all rewritings concerning $\dmd_i$ look the
same for all $i$, we can restrict ourselves to $m=2$ and check the property by
hand or by computer. Since it holds in all that cases, it holds in general.

\subsubsection{Conclusion}

Piecing all the results of this Section together, we get

\begin{theorem}
The free $m$-Dias algebra $Di^{(m)}$ has as series of dimensions the
series of binomial numbers $\binom{n+m-1}{n-1}$.
\end{theorem}

Concerning the operad associated with $m$-Dias algebras, the presentation
given here is quadratic and confluent. So by a known
result~\cite{Dotsenko,Hoffbeck}, this implies

\begin{theorem}
The $m$-Dias operad and the $m$-dendriform operad are both Koszul operads and
dual to each other.
\end{theorem}

One can easily check the first consequence of this fact since the Poincar\'e
series of the $m$-dendriform algebras and their duals are indeed inverses of
each other for the composition of functions: the free $m$-dendriform algebra
satisfies $g(-x) = -x (1+g)^{m+1}$ whereas its dual satisfies
$h(-x) = \frac{-x}{(1+x)^{m+1}}$.


\end{document}